\documentclass{article}
\usepackage{graphicx}
\usepackage{amsfonts}
\usepackage{amsthm}
\usepackage{mathtools}
\usepackage{mathrsfs}
\usepackage{textcomp}
\usepackage{caption}
\usepackage{setspace}
\usepackage{bm}
\usepackage{authblk}
\usepackage{float}
\restylefloat{table}

\newtheorem{definition}{Definition}

\newcommand\extrafootertext[1]{%
    \bgroup
    \renewcommand\thefootnote{\fnsymbol{footnote}}%
    \renewcommand\thempfootnote{\fnsymbol{mpfootnote}}%
    \footnotetext[0]{#1}%
    \egroup
}

\begin{document}
\title{\bf Change of Basis between\\Classical Orthogonal Polynomials}
\markright{Change of Basis}
\author{D.A.\ Wolfram}
\affil{\small College of Engineering \& Computer Science\\The Australian National University, Canberra, ACT 0200\\
\medskip
{\rm David.Wolfram@anu.edu.au}
}
\date{}

\maketitle
\maketitle

\begin{abstract}
Classical orthogonal polynomials have widespread applications including in numerical integration, solving differential equations, and interpolation. Changing basis between classical orthogonal polynomials can affect the convergence, accuracy, and stability of solutions. 

We provide a general method for changing basis between any pair of classical orthogonal polynomials by using algebraic expressions called coefficient functions that evaluate to connection coefficients. The method builds directly on previous work on the change of basis groupoid.  The scope has fifteen kinds of classical orthogonal polynomials including the classes of Jacobi, Gegenbauer and generalized  Laguerre polynomials. 

The method involves the mappings to and from the monomials for these polynomial bases. Sixteen coefficient functions appear to be new  for polynomials that do not have definite parity. We derive the remainder from seven sources in the literature. We give a complete summary of thirty coefficient functions. 

This enables change of basis to be defined algebraically and uniformly between any pair of classical orthogonal polynomial bases by using a  vector dot product to compose two coefficient functions to give a third.  

A definition of Jacobi polynomials uses a basis of shifted monomials. We find a key new mapping between the monomials and Jacobi polynomials  by using a general mapping between shifted monomials and monomials. It yields new mappings for the Chebyshev polynomials of the third and fourth kinds and their shifted versions. 
\end{abstract}

\extrafootertext{MSC: Primary 33C45; Secondary 15A03, 20N02}


\section{Introduction and Related Work}

Change of basis of a finite vector space has numerous significant  applications in scientific computing and engineering.
Applications of change of basis with classical orthogonal polynomials include improving properties of spectral methods for solving differential equations numerically. These properties include better convergence~\cite{boyd2,hale1}, lower computational complexity~\cite{fortunato,olver}, and better numerical stability~\cite{olver}. 

Previous related work gives recurrence relations for finding connection coefficients for orthogonal polynomials, e.g.,~\cite{roncha}, for Jacobi polynomials~\cite{lewanowicz}, or by other methods~\cite{bella,chaggara} and for special cases~\cite{askey2,boyd2}.

This work derives coefficient functions that evaluate to connection coefficients for mappings between the classical orthogonal polynomials and the monomials. We write the coefficient functions as algebraic expressions. From the framework of the change of basis groupoid~\cite{wolfram}, a vector dot product equation enables us to compose coefficient functions in order to derive a coefficient function for any pair of classical orthogonal polynomials.

\section{Classifications of Orthogonal Polynomials}
Thirteen kinds of  orthogonal polynomials are classified by Koornwinder et al.~\cite[\S 18.3]{dlmf} as classical orthogonal polynomials:
\begin{itemize}
\item Jacobi, $P_n^{(\alpha, \beta)}(x)$
\item Gegenbauer, $C_n^{(\lambda)}(x)$
\item Chebyshev of the first to fourth kinds: $T_n(x), U_n(x), V_n(x), W_n(x)$
\item Shifted Chebyshev of the first and second kinds: $T_n^{\ast}(x), U_n^{\ast}(x)$
\item Legendre and Shifted Legendre, $P_n(x), P_n^{\ast}(x)$
\item Generalised Laguerre, $L_n^{(\alpha)}(x)$
\item Physicists' and Probabilists' Hermite: $H_n(x)$, and $He_n(x)$.
\end{itemize}

The Jacobi, Gegenbauer and generalized Laguerre polynomials depend on the respective parameters, $\alpha$ and $\beta$, $\lambda$, and $\alpha$. 

Andrews and Askey~\cite{askey} give a more general classification of the classical orthogonal polynomials and stated ``there are a number of different places to put the boundary for the classical polynomials''. 

We also include the shifted Chebyshev polynomials of the third and fourth kinds here, $V_n^{\ast}(x)$ and $W_n^{\ast}(x)$, because the classification is not fixed.

\section{Related Classical Orthogonal Polynomials}
\label{relatedpoly}

We shall use some of the properties described here that relate different kinds of classical orthogonal polynomials. Specifically, Chebyshev polynomials of the third and fourth kinds and Jacobi polynomials, and Chebyshev polynomials of the second kind and Legendre polynomials as special cases of Gegenbauer polynomials.

The Jacobi and Gegenbauer polynomials are related to some of the other classical orthogonal polynomials. The Gegenbauer polynomials are a special case of the Jacobi polynomials where $\alpha = \beta$, $\alpha = \lambda - \frac12$ and $\alpha > -1$ or $\lambda > -\frac12$. From Szeg\"{o}~\cite[\S 4.7]{szego},
\begin{align}
C_n^{(\lambda)}(x) =& \frac{\Gamma(\alpha+1)}{2 \alpha + 1}\frac{\Gamma(\Gamma(n+ 2 \alpha + 1)}{n + \alpha + 1} P_n^{(\alpha, \alpha)}(x)\\
=& \frac{\Gamma(\lambda + \frac12)}{\Gamma(2\lambda)}\frac{\Gamma(n + 2\lambda)}{\Gamma(n + \lambda + \frac12)} P_n^{(\lambda - \frac12, \lambda - \frac12)}(x)
\end{align}

Chebyshev polynomials of the first kind are related to the Gegenbauer polynomials, e.g., Arfken~\cite{arfken}:

\begin{align}
T_0(x) =& C_0^{(\alpha)}(x) \\
T_n(x) =& \frac{n}2 \lim_{\alpha \rightarrow 0} \frac{C_n^{(\alpha)}(x)}{\alpha} \mbox{ where } n > 0. \label{CT2}
\end{align}

The Chebyshev polynomials of the second kind are a special case of the Gegenbauer polynomials with $U_n(x) = C_n^{(1)}(x)$.  
Chebyshev polynomials of the third and fourth kinds, $V_n(x)$ and $W_n(x)$, are  special cases of the Jacobi polynomials, e.g.,~\cite[\S 1.2.3]{handscomb}.

The Legendre polynomials are also a special case of Gegenbauer polynomials with $P_n(x) = C_n^{(\frac12)}(x)$.

\section{Definitions and Background Equations}

We consider a method of defining algebraically a change of basis mapping between finite bases of classical orthogonal polynomials or the monomials. The method uses the matrix  dot product of change of basis matrices.

From the change of basis groupoid~\cite{wolfram}, given two change of basis matrices $M_{ts} M_{sv}$, we have
\begin{equation} \label {MM}
M_{tv} = M_{ts} M_{sv}.
\end{equation}

The elements of $M_{tv}$ are formed from the vector dot products of the rows of $M_{ts}$ and the columns of $M_{sv}$. 

All classical orthogonal polynomials and the monomials satisfy the following condition. It leads to an optimization of the vector dot products.
\begin{definition} \label{triangular-basis}
Let $V$ be a vector space of finite dimension $n > 0$ and $B_1$ and $B_2$ be bases of $V$.

The basis $B_2$ is a triangular basis with respect to $B_1$ if and only if there is a permutation $\{v_1, \ldots , v_n \}$ of the coordinate vectors of $B_2$ with respect to $B_1$ such that the  $n \times n$ matrix $\left[ v_1^T, \ldots , v_n^T \right]$ is an upper triangular matrix.
\end{definition}

\begin{definition}\label{coefffunction}
Suppose that $s$ and $t$ are triangular bases of a vector space $V$.
The mapping $\mathcal{T}: s \rightarrow t$ satisfies

\begin{equation} \label{mapF}
s_n= \sum_{k=0}^{n} \alpha(n, k)  t_{n - k}
\end{equation}where $\alpha(n, k) \in \mathbb{R}$, i.e., each basis vector of $s$ is a unique linear combination of the basis vectors of $t$.
The function $\alpha$ is called a coefficient function and it evaluates to a connection coefficient.
\end{definition}

If $\alpha(n,k)$ is the coefficient function in equation (\ref{mapF}), then the elements of the change of basis matrix $M_{ts}$ are defined by
\begin{equation} \label{Malpha}
m_{ij} =   \left\{
\begin{array}{ll} 
\alpha(j, j-i)& \mbox{if $j \geq i$}\\
\\
0 & \mbox{if  $j < i$.}\\
\end{array}
\right.
\end{equation} where $0 \leq i, j \leq n$.

Suppose that $s$, $t$ and $v$ are triangular bases,   $\alpha_1(n, k)$  is the coefficient function for a change of basis matrix $M_{ts}$ and $\alpha_2(n,k)$ is a coefficient function for $M_{sv}$ following equation (\ref{Malpha}). Then
\begin{equation} \label{alpha-new2}
\alpha_3(n, k)= \sum_{v=0}^{k} \alpha_1(n - v, k-v) \alpha_2(n , v)
\end{equation} where $0 \leq k \leq n$, is the coefficient function for $M_{tv}$.

The basis $s$ is called the exchange basis~\cite{wolfram}. We shall use the monomials as the exchange basis.

\subsection{Parity and Bases}

Six of the classical orthogonal polynomials have definite parity, and nine do not: the Jacobi polynomials, the generalized Laguerre polynomials, Chebyshev polynomials of the third and fourth kinds, and shifted orthogonal polynomials:  $T^{\ast}$, $U^{\ast}$,  $V^{\ast}$, $W^{\ast}$ and $P^\ast$.

For example, 
 the Jacobi Polynomials $P_n^{(\alpha, \beta)}(x)$ can be expressed in terms of the basis $\{ 1, {(\frac{x-1}2)}, {(\frac{x-1}2)}^2, \ldots, {(\frac{x-1}2)}^n\}$, e.g.,~\cite[equation 18.5.7]{dlmf}. 
The generalised Laguerre Polynomials $L_n^{(\alpha)}(x)$ can be expressed in terms of the basis $\{1, x, x^2, \ldots, x^n\}$.

The Chebyshev Polynomials of the second kind have definite parity and they can be defined by
\[
U_n(x) = 2^n \sum_{k = 0}^{\lfloor \frac{n}{2} \rfloor} (-1)^k  2^{- 2k}{{n-k} \choose k} x^{n - 2k}. 
\] where $n \geq 0$.
When $n$ is even, $U_n(x)$ can be expressed in terms of the basis $\{1, x^2, x^4, \ldots, x^n \}$. When $n$ is odd, $U_n(x)$ can be expressed in terms of $\{x, x^3, x^5, \ldots, x^n\}$.

Equation~(\ref{alpha-new2}) can be optimized when the bases $t$ and $v$ have definite parity and we are concerned with finding the coefficient function for $M_{tv}$ with basis vectors that all have either even parity or odd parity. In the example above, we have a coefficient function $\beta(n, k)$ for non-zero coefficients:
\[
\beta(n,k) =  (-1)^k  2^{n- 2k}{{n-k} \choose k}
\] where $ 0 \leq k  \leq \lfloor \frac{n}{2} \rfloor$.

Given $\beta_1(n,k)$ and $\beta_2(n,k)$ corresponding to the change of basis matrices $M_{ts}$ and $M_{sv}$ respectively, we have
\begin{equation} \label{beta-new}
\beta_3(n, k)= \sum_{v=0}^{k} \beta_1(n - 2v, k-v) \beta_2(n , v)
\end{equation} where $0 \leq k \leq \lfloor \frac{n}2 \rfloor$.   This optimization excludes  terms that are zero~\cite[\S 5.1]{wolfram}.

From equation~(\ref{MM}), $M_{tv} = M_{ts} M_{sv}$, parity is not relevant if at least one of $t$ or $v$ does not have definite parity. In this case,  the following equation can be used when  $\beta$ is known, but $\alpha$ is not.
\begin{equation} \label{beta-alpha}
\alpha(n, k)  =  \left\{
\begin{array}{ll} 
\beta(n, \frac{k}2)& \mbox{if $k$  is even}\\
\\
0 & \mbox{if  $k$ is odd.}\\
\end{array}
\right.
\end{equation}

The upper triangle of the associated change of basis matrix has zero and non-zero alternating elements.

\section{Classical Orthogonal Polynomials and the\\ Monomials}

We summarize here results of mappings between the classical orthogonal polynomials and the monomials. With the monomials as the exchange basis, we can find mappings between any pair of classical orthogonal polynomials.

In the first  table below, the domain is in the first row of a table and the range is the monomials. In the second table, the range is in the first row and the domain is the monomials. Where a citation is given, a related formula has been given in the literature.   The boldface  equation numbers in the tables  refer to equations derived below.  

\newpage
\begin{table}[h]

\centering
\resizebox{\textwidth}{!}{%
\begin{tabular}{|l|l|l|l|l|l|l|l|l|l|l|l|l|l|l|l|}\hline
&$P^{(\alpha, \beta)}$& $C^{(\lambda)}$& $T$&  $U$ & $V$ & $W$ & $T^{\ast}$ & $U^{\ast}$ &  $V^{\ast}$ & $W^{\ast}$ & $P$ & $P^{\ast}$  & $L^{(\alpha)}$ & $H$  & $H_e$\\ \hline
$M$ &{\bf (\ref{Pz1})}& \cite{as}& \cite{as,handscomb}& \cite{as} &\cite{dewi} \& {\bf (\ref{Vx}) }&\cite{dewi} \& {\bf (\ref{Wx})} &{\bf (\ref{MTs})}&{\bf (\ref{MUs})} &{\bf (\ref{MVs})} & {\bf (\ref{MWs})} & \cite{as}&{\bf (\ref{MPs})} & \cite{as} & \cite{as} &  \cite{as} \\ \hline
\end{tabular}}
\caption{Mappings to the Monomials} \label{M1}
\end{table}

\begin{table}[h]
\centering
\resizebox{\textwidth}{!}{%
\begin{tabular}{|l|l|l|l|l|l|l|l|l|l|l|l|l|l|l|l|}\hline
&$P^{(\alpha, \beta)}$& $C^{(\lambda)}$& $T$&  $U$ & $V$ & $W$ & $T^{\ast}$ & $U^{\ast}$ &  $V^{\ast}$ & $W^{\ast}$ & $P$ & $P^{\ast}$  & $L^{(\alpha)}$ & $H$  & $H_e$\\ \hline
$M$ &{\bf (\ref{zP})} & \cite{as,dlmf}& \cite{handscomb,tao}& \cite{dlmf} &{\bf (\ref{xV})}& {\bf (\ref{xW})}&{\bf (\ref{Ts1M}) \& (\ref{Ts2M})} &{\bf (\ref{UsM})} & {\bf (\ref{VsM})}& {\bf (\ref{WsM})}& \cite{dlmf,mathworld1}&  {\bf (\ref{PsM})} & \cite{dlmf} & \cite{dlmf,wiki-hermite} &  \cite{wiki-hermite} \\ \hline
\end{tabular}}
\caption{Mappings from the Monomials} \label{M2}
\end{table}

We are concerned with defining algebraic expressions for the elements in change of basis matrices.
From Tables \ref{M1} and \ref{M2} above, we can use equation~(\ref{alpha-new2}) to find the coefficient function for any pair of bases that are classical orthogonal polynomials.

\subsection{Shifted Monomials and Polynomials}

The  shifted classical orthogonal polynomials and the Jacobi polynomials  are sometimes defined, e.g., \cite[Table 18.3.1]{dlmf} and \cite[\S 1.3.1]{handscomb}, by using a basis of shifted monomials. To express them in terms of the monomials involves changes of bases between shifted monomials and monomials, which we now discuss.

The intervals of orthogonality of the shifted classical  orthogonal polynomials are $(0, 1)$ instead of $(-1, 1)$. They can be defined by applying the polynomial function for the unshifted polynomials to $2 x - 1$ instead of $x$.

Similarly, the interval of orthogonality of the  Jacobi polynomials is $(-1, 1)$~\cite[Table 18.3.1]{dlmf}. They can be viewed as shifted polynomials where the polynomial function has been applied to $\frac{x-1}2$ instead of  $x$, so that the interval before the application was $(-1, 0)$.

Generally, applying a classical orthogonal polynomial function to $c x + d$ instead of $x$ results in a shift of the interval of orthogonality from $(a, b)$ to  $(\frac{a-d}c, \frac{b-d}c)$.

To find the change of basis mappings between the shifted polynomials and the monomials, we first find the change of basis mappings between the monomials and the basis $\{1, cx + d, (cx+d)^2, \ldots, (cx+d)^n\}$.
From the Binomial Theorem, the mapping to the monomials is given by the coefficient function
\begin{equation} \label{Shiftm}
\alpha(n, k, c, d) = {{n} \choose {k}} c^{n-k} {d}^k.
\end{equation}

For example, when $n=3$, $c = \frac13$ and $d = -\frac23$ the change of basis matrix is
\begin{gather*}
\begin{bmatrix}
1 & -\frac23 & \phantom{-}\frac49 & -\frac{8}{27}\\
\\
0 &\phantom{-}\frac13 & -\frac49 & \phantom{-}\frac49\\
\\
0 & \phantom{-}0 & \phantom{-}\frac{1}9 & -\frac{2}{9}\\
\\
0 & \phantom{-}0 & \phantom{-}0 & \phantom{-}\frac{1}{27}
\end{bmatrix}
\end{gather*} and from the fourth column 
\begin{align*}
\left(\frac{x-2}3 \right)^3 
=& \frac{1}{27} x ^3 - \frac29 x^2 + \frac49 x - \frac{8}{27}.
\end{align*}

The inverse mapping from the monomials is given by the coefficient function 
\begin{equation} \label{mShift}
\alpha(n, k, \frac1{c}, -\frac{d}{c}) = {{n} \choose {k}} c^{-n}(-d)^{k}.
\end{equation}

Similarly, when $n=3$,  $c = \frac13$ and $d = -\frac23$, the change of basis matrix is
\begin{gather*}
\begin{bmatrix}
1 & 2 & 4 & 8\\
\\
0 & 3 & 12 & 36\\
\\
0 & 0 & 9 & 54\\
\\
0 & 0 & 0 & 27
\end{bmatrix}
\end{gather*} and from the fourth column
\[
x^3 = 27 \left(\frac{x-2}3\right)^3  + 54 \left(\frac{x-2}3\right)^2 + 36 \left(\frac{x-2}3\right) + 8.
\]
Equations~(\ref{Shiftm}) and~(\ref{mShift}) can be used to give change of basis matrices for the mappings between the monomials and the four kinds of shifted Chebyshev polynomials, and the shifted Legendre polynomials.

For example, with the shifted Chebyshev polynomials of the second kind, we have $U_n^{\ast}(x) = U_n(2x -1)$, so that $c = 2$ and $d = -1$. 
The mapping from the basis of shifted Chebyshev polynomials of the second kind to the monomials follows from equations~(\ref{Shiftm}), Abramowitz and Stegun~\cite[equation 22.3.7]{as} and equation (\ref{alpha-new2}). Specifically, from equation~(\ref{Shiftm}),
\[
\alpha_1(n, k) =  
{n \choose k} 2^{n -k}  {(-1)}^{k}.
\]

From equations 22.3.7 and (\ref{beta-alpha}),
\[
\alpha_2(n, k)=  \left\{
\begin{array}{ll} 
{{n-\frac{k}2} \choose \frac{k}2} 2^{n -k}   {(-1)}^{\frac{k}2} & \mbox{if $k$ is even}\\
\\
0 & \mbox{if  $k$ is odd.}\\
\end{array}
\right.
\]

Hence from equation~(\ref{alpha-new2}), when $n=3$ the change of basis matrix is
\begin{gather*}
\begin{bmatrix}
1 & -2 & \phantom{-}3 & -4\\
\\
0 & \phantom{-}4 & -16 & \phantom{-}40\\
\\
0 & \phantom{-}0 & \phantom{-}16 & -96\\
\\
0& \phantom{-}0 & \phantom{-}0& \phantom{-}64
\end{bmatrix}
\end{gather*} and from the fourth column
\[
U_3(2 x - 1) = 64 x^3 - 96 x^2 + 40 x -4.
\]

The inverse mapping from the monomials to the shifted Chebyshev or shifted Legendre polynomials uses equation~(\ref{mShift}), the mapping from monomials to the unshifted polynomials, and equation~(\ref{alpha-new2}).  

For example, for the shifted Chebyshev polynomials of the first kind, we use the mapping from the monomials to the Chebyshev Polynomials of the first kind~\cite{handscomb,tao} and equation~(\ref{beta-alpha}) to give
\[
\alpha_1(n, k) =  \left\{
\begin{array}{ll} 
{{n} \choose \frac{k}2} 2^{1-n}   & \mbox{if $k$ is even and $0 \leq k < n$}\\
\\
{{n} \choose \frac{k}2} 2^{-n}   & \mbox{if $k$ is even and $k= n$}\\
\\
0 & \mbox{if  $k$ is odd.}\\
\end{array}
\right.
\] 
With $c=2$ and $d = -1$, equation~(\ref{mShift}) gives
\[
\alpha_2(n, k) = 2^{-n} {n \choose k}.
\]

The change of basis matrix is
\begin{gather*}
\begin{bmatrix}
1 & \frac12 & \frac38 & \frac5{16}\\
\\
0 & \frac12 & \frac12 & \frac{15}{32}\\
\\
0 & 0 & \frac18 & \frac{3}{16}\\
\\
0& 0 & 0& \frac1{32}
\end{bmatrix}
\end{gather*} whose elements are 
\[
\alpha_3(j, j-i) = \sum_{v=0}^{j-i} \alpha_1(j - v, j-i-v) \alpha_2(j, v)
\] where $0 \leq i \leq j \leq 3$ and $0$ otherwise. From the third and fourth columns, we have
\begin{align*}
x^2 =& \frac{1}8 T_2(2x -1) + \frac12 T_1(2 x -1) + \frac38\\
x^3 =& \frac{1}{32}T_3(2x -1) + \frac{3}{16}T_2(2x -1) + \frac{15}{32} T_1(2x -1)  + \frac{5}{16}.
\end{align*}

\subsection{Jacobi Polynomials}
A definition of the Jacobi polynomials~\cite[equation 18.5.7]{dlmf}  is
\begin{equation*}
P_n^{(\alpha, \beta)}(z) = \sum_{l=0}^n\frac{(n+\alpha+\beta+1)_l (\alpha+l+1)_{n-l}}{l! (n-l)!} {\left(\frac{z-1}2\right)}^l
\end{equation*}
which uses the Pochhammer symbol or rising factorial.
This can be defined as a polynomial in the basis $\{1, \frac{z-1}2, {\left(\frac{z-1}2\right)}^2, \ldots \}$. This basis is a triangular basis following Definition~(\ref{triangular-basis}). We can change basis to the monomials.

From equation (\ref{Shiftm}), the change of basis matrices from $\{1, \frac{z-1}2, {\left(\frac{z-1}2\right)}^2, \ldots \}$ to monomials has coefficient function
\begin{equation}
\alpha(n, k) = {{n} \choose {k}} 2^{-n}(-1)^{k}
\end{equation} where $c = \frac12$ and $d = -\frac12$.

For example, when $n = 4$, the change of basis matrix  is
\begin{gather*}
\begin{bmatrix}
1 & -\frac12 & \frac14 & \frac18\\
\\
0 & \frac12 & -\frac12 & \frac38\\
\\
0 & 0 & \frac14 & -\frac38\\
\\
0& 0 & 0& \frac18
\end{bmatrix}
\end{gather*}

The coordinates of $P_3^{(\alpha, \beta)}(z)$ in the domain basis written as a vector are
\begin{gather*}
\begin{bmatrix}
\frac16 (1+ \alpha)(2+ \alpha)(3+ \alpha)\\
\\
\frac12 (2 + \alpha)(3+\alpha)(4 + \alpha+ \beta) \\
\\
\frac12 (3 + \alpha)(4 + \alpha + \beta)(5 + \alpha + \beta)\\
\\
\frac16 (4 + \alpha + \beta) (5 + \alpha + \beta) (6 + \alpha + \beta)
\end{bmatrix}
\end{gather*}

The product of the change of basis matrix above with this vector, gives the coordinates of the Jacobi polynomial in the basis of the monomials expressed as a vector.

For example, the expression below is the third element, and it is the coefficient of the term in $z^2$ of $P_3^{(\alpha, \beta)}(z)$:
\[
\frac18  (3 + \alpha)(4 + \alpha + \beta)(5 + \alpha + \beta) -\frac{1}{16} (4 + \alpha + \beta) (5 + \alpha + \beta) (6 + \alpha + \beta)
\] which equals
\[
\frac{5a}4 + \frac{9 \alpha^2}{16} + \frac{\alpha^3}{16} - \frac{5\beta}4 + \frac{\alpha^2 \beta}{16} - \frac{9\beta^2}{16} - \frac{\alpha\beta^2}{16} - \frac{\beta^3}{16}.
\]

It follows that 
\begin{align} \label{Pz}
P_n^{(\alpha, \beta)}(z) =& 
\sum_{m=0}^n \left( \sum_{l=0}^{n-m} \frac{ (-1)^{n-m}}{2^{n-l}}   \frac{(1 + \alpha+\beta+n)_{n-l}(-\alpha - n)_{l}}{(n-m-l)! \; m! \; l!} \right) z^m
\end{align}

or equivalently
\begin{align} 
P_n^{(\alpha, \beta)}(z) =& 
\frac{\Gamma(\alpha+n+1)}{n!\Gamma(\alpha+\beta+n+1)} \nonumber\\
&\sum_{m=0}^n \left( \sum_{l=0}^{n-m}(-2)^{l - n}  (-1)^{m} {{n-l} \choose m} {n \choose {n-l}}  \frac{\Gamma(\alpha+\beta+2n - l +1)}{\Gamma(\alpha + n-l + 1)} \right) z^m
\end{align}

Therefore, we have,
\begin{align}
\alpha(n,k) =& \frac{(-1)^{n-k} \Gamma(\alpha+n+1)}{n!\Gamma(\alpha+\beta+n+1)} \nonumber\\
&\sum_{l=0}^{k}(-2)^{l - n}  {{n-l} \choose {n-k}} {n \choose {n-l}}  \frac{\Gamma(\alpha+\beta+2n - l +1)}{\Gamma(\alpha + n-l + 1)}. \label{Pz1}
\end{align}

For the inverse mapping, equation 18.18.15 of Koornwinder et al.~\cite{dlmf} is
\[
\left( \frac{1+x}2 \right)^n = \frac{(\beta+1)_n}{(\alpha + \beta + 2)_n} \sum_{l = 0}^n \frac{\alpha+\beta + 2l + 1}{\alpha + \beta + 1}
\frac{(\alpha + \beta+1)_l (n - l + 1)_l}{(\beta+1)_l (n + \alpha + \beta + 2)_l} P_l^{(\alpha, \beta)}(x).
\] 

This defines a change of basis from  $\{1, \frac{1+x}2, {\left(\frac{1+x}2\right)}^2, \ldots \}$ to the Jacobi polynomials. 
From equation (\ref{mShift}), the change of basis matrices from  the  monomials  to $\{1, \frac{1+x}2, {\left(\frac{1+x}2\right)}^2, \ldots \}$ has coefficient function
\begin{equation} \label{mJacobi}
\alpha(n, k) = {{n} \choose {k}} 2^{n-k} (-1)^k
\end{equation} where $c = \frac12$ and $d = \frac12$. This gives
\begin{align} 
z^n =&  \sum_{m=0}^{n} \frac{\alpha+\beta +2m + 1}{\alpha + \beta + 1}  \left( 
\sum_{l = 0}^{n -m} 2^{n-l} {n \choose {n-l}} {(-1)}^l  \frac{(\beta+1)_{n-l}}{(\alpha + \beta + 2)_{n-l}}\right. \nonumber \\
&\left.  \frac{(\alpha + \beta+1)_{m} (n - l - m + 1)_{m}}{(\beta+1)_{m} (n-l+ \alpha + \beta + 2)_{m}}
\right) 
P_{m}^{(\alpha, \beta)}(z).
\end{align}

We have
\begin{align} \label{zP}
\alpha(n,k)  =&  \ (\alpha+\beta +2(n-k) + 1)(\alpha + \beta + 2)_{n-k-1}  \nonumber \\
 &\sum_{l = 0}^{k} 2^{n-l} {n \choose {n-l}} {(-1)}^l  \frac{(\beta+1+n-k)_{k-l}}{(\alpha + \beta + 2)_{n-l}} 
\frac{ (k - l  + 1)_{n-k}}{(n-l+ \alpha + \beta + 2)_{n-k}}.
\end{align}

\subsection{Chebyshev Polynomials of the Third and Fourth Kinds}

These classical orthogonal polynomials can be respectively defined from the Jacobi polynomials~\cite{handscomb}
\begin{align}
V_n(x) =& \frac{2^{2n}}{{{2n} \choose n}} P_n^{(-\frac12, \frac12)}(x) \label{VP} \\ 
W_n(x) =& \frac{2^{2n}}{{{2n} \choose n}} P_n^{(\frac12, -\frac12)}(x).  \label{WP}
\end{align}  $V_n(x)$ and $W_n(x)$ have terms in all powers of $x$ up to $n$ when expressed in terms of monomials.

From equation~(\ref{Pz}), we have 
\begin{align} 
V_n(x) =& 
\frac{2^{n}}{{{2n} \choose n}} \sum_{m=0}^n  (-1)^{n-m} \left( \sum_{l=0}^{n-m} 2^l    \frac{(1 +n)_{n-l}(\frac{1}{2} - n)_{l}}{(n-l- m)! \; m! \; l!} \right) x^m 
\end{align} so that
\begin{align} \label{Vx}
\alpha(n,k) =& \frac{2^{n}}{{{2n} \choose n}}  \frac{(-1)^{k}}{ (n-k)!} \sum_{l=0}^{k} 2^l    \frac{(1 +n)_{n-l}(\frac{1}{2} - n)_{l}}{(k-l)! \; l!}. 
\end{align}

Similarly,
\begin{align}
W_n(x) =& 
\frac{2^{n}}{{{2n} \choose n}} \sum_{m=0}^n  (-1)^{n-m} \left( \sum_{l=0}^{n-m} 2^l    \frac{(1 +n)_{n-l}(-\frac{1}{2} - n)_{l}}{(n-l-m)! \; m! \; l!} \right) x^m 
\end{align}  so that
\begin{align}\label{Wx}
\alpha(n,k) =&  
\frac{2^{n}}{{{2n} \choose n}}  \frac{(-1)^{k}}{(n-k)!}  \sum_{l=0}^{k} 2^l    \frac{(1 +n)_{n-l}(-\frac{1}{2} - n)_{l}}{(k-l)! \;  l!}.
\end{align}

Explicit formulas for these polynomials were also given by another approach by Dewi, Utama and Animah~\cite{dewi}. They are
\begin{align*} 
V_n(x) =& \sum_{k = \lceil \frac{n}2 \rceil}^n {k \choose {n-k}} 2^{2k -n -1}(-1)^{n-k}  x^{2k -n -1} \left(2x - \frac{2k-n}k\right)\\
W_n(x) =& \sum_{k = \lceil \frac{n}2 \rceil}^n {k \choose {n-k}} 2^{2k -n -1}(-1)^{n-k}  x^{2k -n -1} \left(2x + \frac{2k-n}k\right)
\end{align*} where $n \geq 2$. These formulas also apply when $n=1$.

In the converse direction, from equation~(\ref{zP}), let
\[
 f(m,n) = \sum_{l = 0}^{n -m} 2^{n-l} {n \choose {n-l}} {(-1)}^l  \frac{(\frac32)_{n-l}}{(2)_{n-l}} \frac{(1)_{m} (n - l - m + 1)_{m}}{(\frac32)_{m} (n-l+ 2)_{m}}
 \] 
so that from equation~(\ref{VP}), we obtain
\begin{equation}
x^n =   \sum_{m=0}^{n} (2m + 1) \frac{{{2m} \choose m}}{2^{2m}} f(m,n) V_m(x).
\end{equation}

This can be simplified by using 
\[
(x)_n = \frac{\Gamma(x + n)}{\Gamma(x)} \:\:\mbox{ and }\:\: \Gamma(\frac12 + n) = \frac{(2n - 1)!!}{2^n}\sqrt{\pi}
\] to give

\begin{align}
x^n =&   \sum_{m=0}^{n} \left(\sum_{l = 0}^{n -m} \frac{ (-1)^l  \;(2(n-l) +1)!! \; n!}{(n-l- m)! \; (n-l+m+1)! \; l!} \right) V_m(x).
\end{align}

The coefficient function for the change of basis from monomials to Chebyshev polynomials of the third kind is therefore
\begin{equation}  \label{xV}
\alpha(n, k) = n! \sum_{l = 0}^{k} \frac{ (-1)^l  \;(2(n-l) +1)!!}{(k-l)! \; (2n- l -k  +1)! \; l!}.
\end{equation}

Similarly, we obtain
\begin{align} 
x^n =&   \sum_{m=0}^{n} \left(\sum_{l = 0}^{n -m}  (2m+1)\frac{(-1)^l \; (2(n-l)-1)!! \; n!}{(n-l-m)! \; (n-l+m+1)! \;  l!} \right) W_m(x).
\end{align}

The coefficient function for the change of basis matrix from $W_n(x)$ to the monomials is
\begin{equation} \label{xW}
\alpha(n, k) =  
(2(n-k)+1) n! \sum_{l = 0}^{k}  \frac{(-1)^l \; (2(n-l)-1)!! }{(k-l)! \; (2n-l -k +1)! \;  l!}.
\end{equation} 

For example, when $n=4$, we have

\begin{gather*}
\begin{bmatrix}
1 & -\frac12 & \phantom{-}\frac12 & -\frac38 & \phantom{-}\frac38\\
\\
0 & \phantom{-}\frac12& -\frac14 & \phantom{-}\frac38 & -\frac14\\
\\
0 & 0 & \phantom{-}\frac14& -\frac18 & \phantom{-}\frac14\\
\\
0& 0 & 0& \phantom{-}\frac18& -\frac1{16}\\
\\
0&0&0&0&\phantom{-}\frac{1}{16}
\end{bmatrix}
\end{gather*}
and from the fifth column, for example,
\begin{align*}
x^4 =& \frac1{16} W_4(x) -\frac1{16}W_3(x)  + \frac14 W_2(x) - \frac14 W_1(x) + \frac38 W_0(x).
\end{align*}

\section{Shifted Classical Orthogonal Polynomials}

The shifted classical orthogonal polynomials include the shifted Chebyshev polynomials of the first kind, $T^{\ast}_n(x)$, the fourth kind, $W^{\ast}_n(x)$, and the shifted Legendre polynomials,  $P^{\ast}_n(x)$. All of them can be defined by applying the unshifted polynomial function to the argument $2 x -1$ instead of $x$, e.g., $T_n^{\ast}(x) = T_n(2 x -1)$. They each have $(0, 1)$ as their domain of orthogonality~\cite[Table 18.3.1]{dlmf}. We note that the method here is general and can be used for any linear shift of the form $cx + d$ and for any triangular polynomial basis.

\subsection{Mappings to the Monomials}
With these definitions, and the definitions of the unshifted polynomials in terms of the monomials, the shifted polynomials can be expressed using the basis $\{1, 2x -1, (2x -1)^2, \ldots, (2x -1)^{n}\}$. To express the shifted polynomials in terms of the monomials, we use equation~(\ref{Shiftm}) with $c = 2$ and $d = -1$, i.e.
\begin{equation} \label{Mshift}
\alpha(n, k) = {n \choose k} (-1)^k 2^{n-k} 
\end{equation} and then apply equation~(\ref{beta-new}) because these polynomials have definite parity.

For shifted Chebyshev polynomials of the first kind we have 
\[
T_n(x) = \frac{n}2 2^n \sum_{k=0}^{\lfloor \frac{n}2 \rfloor} (-1)^k 2^{-2k} \frac{(n-k-1)!}{k!(n-2k)!} x^{n-2k}
\] so that
\[
\beta_2(n,k) = \frac{n}2 2^n (-1)^k 2^{-2k} \frac{(n-k-1)!}{k!(n-2k)!}.
\]

Therefore,
\begin{align}
\alpha(n,k) =& \sum_{v=0}^{\lfloor \frac{k}2 \rfloor}  {{n-2v} \choose {k-2v}} (-1)^{k-2v} 2^{n-k} \frac{n}2 2^n (-1)^{v} 2^{-2v} \frac{(n-v-1)!}{v!(n-2v)!} \nonumber\\
=& \frac{n}2  \sum_{v=0}^{\lfloor \frac{k}2 \rfloor}  {{n-2v} \choose {k-2v}}  (-1)^{k-v} 2^{2(n-v)-k}   \frac{(n-v-1)!}{v!(n-2v)!} \label{MTs}
\end{align} and $\alpha(0,0) = 1$. We have
\[
T_n^{\ast}(x) = \sum_{k=0}^{n} \alpha(n,k) x^{n-k}.
\]

Similarly, for shifted Chebyshev polynomials of the second kind, we obtain
\begin{align}
\alpha(n,k) =& \sum_{v=0}^{\lfloor \frac{k}2 \rfloor}  {{n-2v} \choose {k-2v}} (-1)^{k-2v} 2^{n-k} 2^{n-2v} {{n-v}\choose v} (-1)^v \nonumber\\
=&\sum_{v=0}^{\lfloor \frac{k}2 \rfloor}  {{n-2v} \choose {k-2v}} {{n-v}\choose v} (-1)^{k-v} 2^{2(n-v)-k}.  \label{MUs}
\end{align}

The shifted Chebyshev polynomials of the third and fourth kinds, $V^{\ast}$ and $W^{\ast}$, can be similarly defined by
$V_n^{\ast}(x) = V_n(2 x -1)$  and $W_n^{\ast}(x) = W_n(2x -1)$, e.g.,~\cite[equation (1.27)]{handscomb}

For the mapping from the monomials to $V_n^{\ast}(x)$, we apply equation~(\ref{alpha-new2}) with $\alpha_1$ from equation~(\ref{Mshift}) and $\alpha_2$ from equation~(\ref{Vx}) to give :
\begin{equation} \label{MVs}
\alpha_3(n, k) = (-1)^{k} \frac{2^{2n-k}}{{{2n} \choose n} (n-k)!}  \sum_{v=0}^k
  \sum_{l=0}^{v} 2^l    \frac{(1 +n)_{n-l}(\frac{1}{2} - n)_{l}}{(k-v)! \; (v-l)! \;  l!}.
\end{equation}

For the similar mapping to $W_n^{\ast}(x)$ we use $\alpha_2$ from equation~(\ref{Wx}) instead to obtain:
\begin{equation} \label{MWs}
\alpha_3(n, k) = (-1)^{k} \frac{2^{2n-k}}{{{2n} \choose n} (n-k)!}  \sum_{v=0}^k
  \sum_{l=0}^{v} 2^l    \frac{(1 +n)_{n-l}(-\frac{1}{2} - n)_{l}}{(k-v)! \; (v-l)! \;  l!}.
\end{equation}

For shifted Legendre polynomials, we have
\begin{align}
\alpha(n,k) =& \sum_{v=0}^{\lfloor \frac{k}2 \rfloor}  {{n-2v} \choose {k-2v}} (-1)^{k-2v} 2^{n-k} 2^{-n}  {{2n-2v}\choose n} {n \choose v}(-1)^v \nonumber\\
=&2^{-k} \sum_{v=0}^{\lfloor \frac{k}2 \rfloor}  {{n-2v} \choose {k-2v}} {{2n-2v}\choose n} {n \choose v} (-1)^{k-v}.   \label{MPs}
\end{align} 

\subsection{Mappings from the Monomials}
In the converse direction, we apply the  change of bases from 
\[
\{1, 2x -1, (2x -1)^2, \ldots, (2x -1)^{n}\}
\] to the shifted polynomials $T_n^{\ast}(x)$, $U_n^{\ast}(x)$ and $P_n^\ast(x)$. We first consider the change of basis from the monomials to $\{1, 2x -1, (2x -1)^2, \ldots, (2x -1)^{n}\}$.

This mapping follows from equation (\ref{mShift}) with $c = 2$ and $d = -1$ to give
\begin{equation} \label{invShift}
\alpha_2(n, k) = {n \choose k} 2^{-n}.
\end{equation}
 
The change of basis from monomials to  Chebyshev polynomials of the first kind is given by  the coefficient function
\begin{equation*}
\beta_1(n, k) =  \left\{
\begin{array}{ll}
{n \choose k} 2^{1-n}& \mbox{if $0  \leq 2 k< n$}\\ 
\\
{n \choose k} 2^{-n} &  \mbox{if $2 k  = n$.}
\end{array}
\right.
\end{equation*}

This is equivalent to solutions given by Mason and Handscomb~\cite{handscomb}, section 2.3.1 and a note by Tao~\cite{tao}.

The coefficient function is the same as that of shifted Chebyshev functions of the first kind since the variable $x$ is replaced by $2x -1$. These polynomials have definite parity, so that the optimization in case 3 of Wolfram~\cite[equation (12)]{wolfram} can be applied to give the required coefficient function:
\begin{align}
\alpha_3(n, k)=& \sum^{k}_{\mathclap{\substack{v = 0 \\ k-v \mbox{\small \; even}}}} \alpha_1(n  - v, k-v) \alpha_2(n , v) \nonumber\\
=& \sum^{k}_{\mathclap{\substack{v = 0 \\ k-v \mbox{\small \; even}}}} \beta_1(n  - v,  \frac{k-v}2) {n \choose v} 2^{-n}. \nonumber
\end{align} This gives
\begin{equation}
\alpha_3(n,k) = 2^{1-2n} \sum^{k}_{\mathclap{\substack{v = 0 \\ k-v \mbox{\small \; even}}}} {{n-v} \choose { \frac{k-v}2}} {n \choose v}  2^v\:\:\: \mbox{ if $0  \leq  k< n$} \label{Ts1M}
\end{equation}

\begin{equation}
\alpha_3(n,k) = 2^{-2n} \sum^{k}_{\mathclap{\substack{v = 0 \\ k-v \mbox{\small \; even}}}} {{n-v} \choose { \frac{k-v}2}}    {n \choose v}  2^v \:\:\: \mbox{ if $k  = n$} \label{Ts2M}
\end{equation} and
\[
x^n = \sum_ {k=0}^n \alpha_3(n, k) T_{n-k}^{\ast}(x).
\]

For the  Chebyshev polynomials of the second kind, the mapping from the monomials has the coefficient function
\begin{equation}
\beta_1(n, k) = 2^{-n} (n - 2k + 1) \frac{n!}{(2)_{n-k} k!}.
\end{equation} This follows from Koornwinder et al.~\cite[equation 18.18.17]{dlmf} with $\lambda = 1$. In a similar way to the case for $T_n^{\ast}(x)$ above, we have
\begin{align}
\alpha_3(n, k)=& \sum^{k}_{\mathclap{\substack{v = 0 \\ k-v \mbox{\small \; even}}}} \alpha_1(n  - v, k-v) \alpha_2(n , v) \nonumber\\
=& \sum^{k}_{\mathclap{\substack{v = 0 \\ k-v \mbox{\small \; even}}}} \beta_1(n  - v,  \frac{k-v}2) {n \choose v} 2^{-n} \nonumber \\
=&  2^{-2n}(1 + n -k) n! \sum^{k}_{\mathclap{\substack{v = 0 \\ k-v \mbox{\small \; even}}}}   \frac{2^v}{(2)_{n - \frac{v+k}2} v! (\frac{k-v}2)!}. \label{UsM}
\end{align}

For the mapping from $V^{\ast}_n(x)$ from the monomials, we use $\alpha_1$ from equation~(\ref{xV}) and $\alpha_2$ from equation~(\ref{invShift})  to give:
\begin{equation} \label{VsM}
\alpha_3(n, k) =  2^{-n} n! \sum_{v=0}^k
  \sum_{l = 0}^{k-v} \frac{ (-1)^l  \;(2(n -v-l) +1)!!}{(k-v- l)! \; (2n-v - l - k  +1)! \; v! \; l!}.
\end{equation}

Similarly for $W^{\ast}_n(x)$  we use $\alpha_1$ from equation~(\ref{xW}):
\begin{equation} \label{WsM}
\alpha_3(n, k) =  2^{-n} (2(n-k) + 1) n! \sum_{v=0}^k 
  \sum_{l = 0}^{k-v} \frac{ (-1)^l  \;(2(n -v-l) -1)!!}{(k-v- l)! \; (2n-v - l - k  +1)! \; v! \; l!}.
\end{equation}

For the shifted Legendre polynomials,  we have from Koornwinder et al.~\cite[equation 18.18.17]{dlmf} with $\lambda = \frac12$ that
\begin{equation}
\beta_1(n, k) = 2^{-n} (1 + 2n - 4k) \frac{n!}{(\frac32)_{n-k} k!}
\end{equation} so that
\begin{align}
\alpha_3(n, k)=& \sum^{k}_{\mathclap{\substack{v = 0 \\ k-v \mbox{\small \; even}}}} \alpha_1(n  - v, k-v) \alpha_2(n , v) \nonumber\\
=&  \sum^{k}_{\mathclap{\substack{v = 0 \\ k-v \mbox{\small \; even}}}} \beta_1(n  - v,  \frac{k-v}2) {n \choose v} 2^{-n} \nonumber \\
=&  2^{-2n} (2(n-k) + 1)  n!\sum^{k}_{\mathclap{\substack{v = 0 \\ k-v \mbox{\small \; even}}}}  \frac{2^v}{(\frac32)_{n - \frac{v+k}2} v! (\frac{k-v}2)!}.  \label{PsM}
\end{align}

For example,
\begin{align*}
x^4 =& \alpha_3(4, 0) P_4^{\ast}(x) + \alpha_3(4, 1) P_3^{\ast}(x) + \alpha_3(4, 2) P_2^{\ast}(x) + \alpha_3(4, 3) P_1^{\ast}(x) +\alpha_3(4, 4) P_0^{\ast}(x) \\
=& \frac{1}{70}  P_4(2x -1) + \frac1{10} P_3(2x-1)+ \frac{2}7 P_2(2x-1) +  \frac25 P_1(2x-1)+ \frac15 .
\end{align*}

\section{Examples}

We give two examples of derivations of coefficient functions for change of basis between Jacobi polynomials,  and Physicist's Hermite polynomials to the shifted Legendre polynomials.

The coefficient function for Jacobi polynomials has four extra parameters which are the parameters of the Jacobi polynomials: two for the those in the domain basis and two for the range basis. By swapping these parameters, we find the inverse of the change of basis matrix.

The example of deriving the coefficient function for the change of basis from Physicist's Hermite polynomials to the shifted Legendre polynomials uses equation~(\ref{alpha-new2}) because the shifted Legendre polynomials do not have definite parity.

\subsection{Jacobi Polynomials}

From equation~(\ref{Pz}), we have

\begin{align}
\alpha_2(n, k, \alpha, \beta) =& \sum_{l=0}^{k} \frac{ (-1)^{k}}{2^{n-l}}   \frac{(1 + \alpha+\beta+n)_{n-l}(-\alpha - n)_{l}}{(k-l)! \; (n-k)! \; l!}.
\end{align}

From equation~(\ref{zP}), let $k = n - m$, so that 
\begin{align} \label{alpha-zP}
\alpha_1(n, k, \gamma, \delta) =& \frac{\gamma+\delta +2(n-k) + 1}{\gamma + \delta + 1}  \left( 
\sum_{l = 0}^{k} 2^{n-l} {n \choose {n-l}} {(-1)}^l  \frac{(\delta+1)_{n-l}}{(\gamma + \delta + 2)_{n-l}}\right. \nonumber \\
&\left.  \frac{(\gamma + \delta+1)_{n-k} (k - l  + 1)_{n-k}}{(\delta+1)_{n-k} (n-l+ \gamma + \delta + 2)_{n-k}}
\right).
\end{align}

Hence from equation~(\ref{alpha-new2}), it follows that
\begin{equation}  \label{PP}
\alpha_3(n, k, \alpha, \beta, \gamma, \delta)= \sum_{v=0}^{k} \alpha_1(n - v, k-v,  \gamma, \delta) \alpha_2(n , v, \alpha, \beta).
\end{equation} where $0 \leq k \leq n$.

This equation gives the elements of the upper triangular change of basis matrix for the mapping from the basis  $\{P_n^{(\alpha, \beta)}(z), n \geq 0\}$ to the basis   $\{P_n^{(\gamma, \delta)}(z), n \geq 0\}$.  

Another coefficient function  for Jacobi Polynomials can be derived from equation 18.18.14 of Koornwinder et al.~\cite{dlmf}. These coefficient functions
generalize three solutions for special cases by Askey~\cite{askey2} who used another method: equation (7.32) where $ \alpha = \gamma$, equation (7.33) where $\beta = \delta$, and equation (7.34) where $\alpha = \beta$ and $\gamma = \delta$.

The inverse change of basis matrix has the property that it can be found by interchanging $\alpha$ and $\gamma$ with $\beta$ and $\delta$, respectively, by using $\alpha_3(n, k, \gamma, \delta, \alpha, \beta)$.

For example, let $n = 4$, $\alpha = 1$ and $\beta = 8$, $\gamma = 2$,  $\delta = 7$.  The change of basis matrix from $\{P_n^{(2, 7)}(z), 0 \leq n \leq 4\}$ to the basis  $\{P_n^{(1, 8)}(z), 0 \leq n \leq 4\}$ is

\begin{gather*}
\begin{bmatrix}
1 & 1 & \frac{9}{11} & \frac{15}{22} & \frac{15}{26}\\
\\
0 & 1& \frac{12}{11} & \frac{10}{11} & \frac{10}{13}\\
\\
0 & 0 & 1& \frac76 & \frac{77}{78}\\
\\
0& 0 & 0& 1& \frac{16}{13}\\
\\
0&0&0&0&1
\end{bmatrix}
\end{gather*} and from the fourth column,
\[
P_3^{(2, 7)}(z) = P_3^{(1, 8)}(z) + \frac76 P_2^{(1, 8)}(z) + \frac{10}{11}P_1^{(1, 8)}(z) + \frac{15}{22} P_0^{(1, 8)}(z).
\]

\subsection{Physicist's Hermite Polynomials to Shifted  Legendre Polynomials}

The coefficient function in this example can be found from equation~(\ref{alpha-new2}) because shifted Legendre polynomials do not have definite parity.

The coefficient function for the change of basis from the monomials to shifted Legendre polynomials is 
from equation~(\ref{PsM}) for $\alpha_1(n, k)$.

From Table~\ref{M1} the coefficient function for the change of basis from Physicist's Hermite polynomials to the monomials is derived from Abramowitz and Stegun~\cite[equation 22.3.10]{as}. This is
\begin{equation}
\beta_2(n,k) = (-1)^k 2^{n-2k} \frac{n!}{k!(n-2k)!}.
\end{equation}

From equation~ (\ref{beta-alpha}), 
\begin{equation}
\alpha_2(n, v) =
 \left\{
\begin{array}{ll} 
 (-1)^{\frac{v}2}  2^{n-v} \frac{n!}{(\frac{v}2)!(n-v)!} & \mbox{if $v$  is even}\\
\\
0 & \mbox{if  $v$ is odd.}\\
\end{array}
\right.
\end{equation}

From equation~(\ref{alpha-new2}), we have
\begin{align*}
\alpha_3(n, k)=& \sum_{v=0}^{k} \alpha_1(n - v, k-v) \alpha_2(n , v).
\end{align*}

It follows that
\begin{align*}
H_5(x) =& \sum_{v=0}^5 \alpha_3(5,v) P^{\ast}_{5-v}(x)  \\
=& \frac{8}{63} P^{\ast}_5(x) + \frac87 P^{\ast}_4(x) - \frac{32}9 P^{\ast}_3(x) - \frac{640}{21} P^{\ast}_2(x)  - \frac47 P^{\ast}_1(x)  + \frac{76}3 P^{\ast}_0(x)\\
=& 32 x^5 -160 x^3 + 120 x.
\end{align*}

\section{Conclusion}

Previous work~\cite{wolfram} introduced a change of basis groupoid and provides the framework for this work. We reduced the scope here to fifteen classical orthogonal polynomials including the classes of Jacobi polynomials, Gegenbauer polynomials and generalized Laguerre polynomials. 

We used the monomials  as the exchange basis and gave thirty algebraic expressions called coefficient functions that evaluate to connection coefficients  for the mappings to and from the monomials. The appendix summarizes them.
With this library of coefficient functions, the framework enables us  to derive a coefficient function for the mapping between any pair of classical orthogonal polynomial bases. We do this using equation~(\ref{alpha-new2}), or its optimizations~\cite{wolfram} when at least one of the bases has definite parity.

Of the thirty algebraic expressions, eighteen relate to classical orthogonal polynomials that do not have definite parity. Sixteen of these expressions seem new. 

We used the change of basis between shifted monomials and the monomials. This enabled us to map the basis of Jacobi polynomials to and from the monomials. This mapping led to a similar one for Chebyshev polynomials of the third and fourth kinds by using a definition of them based on the Jacobi polynomials. The technique also produced mappings for five shifted classical orthogonal polynomials: $T^{\ast}$, $U^{\ast}$, $V^{\ast}$, $W^{\ast}$ and $P^{\ast}$. The technique allows for any linear shift $cx + d$, and it applies to any triangular polynomial basis over $\mathbb{C}[x]$.

We then gave examples for changes of bases between Jacobi polynomials, and from Physicist's Hermite polynomials to the shifted Legendre polynomials.

\section*{Acknowledgment}
 I am grateful  to the College of Engineering \& Computer Science at The Australian National University for research support.

\appendix
\section{Coefficient Functions}

We summarize here coefficient functions based on the results above or derived from formulas in the literature. It has three parts: coefficient functions for change of basis for shifted monomials; for the classical orthogonal polynomials to the monomials; and the monomials to the classical orthogonal polynomials.

\subsection{Shifted Monomials}

The mapping to the monomials from the basis $\{1, cx + d, (cx+d)^2, \ldots, (cx+d)^n\}$ is given in the first table below. The second table is for the inverse mapping. 
\begin{table}[h]
\centering
\begin{tabular}{|l|l|l|}\hline
Domain Basis & Equation & $\alpha(n, k, c, d)$\\ \hline
&&\\
$ (cx+d)^n$ & (\ref{Shiftm}) & ${n \choose k} c^{n-k} d^{k}$\\\
&&\\\hline
\end{tabular}
\caption{Mapping to  the Monomials} 
\end{table} 

\begin{table}[h]
\centering
\begin{tabular}{|l|l|l|}\hline
Range Basis & Equation & $\alpha(n, k, c, d)$\\ \hline
&&\\
$ (cx+d)^n$ &(\ref{mShift}) &  ${n \choose k} c^{-n} (-d)^{k}$\\ 
&&\\\hline
\end{tabular}
\caption{Mapping from the Monomials} 
\end{table} 

\newpage
\subsection{Mappings to and from the Monomials for\\Classical Orthogonal Polynomials}

\begin{table}[h]
\centering
\begin{tabular}{|l|l|l|}\hline
Domain Basis & Equation & $\beta(n, k)$\\ \hline
&&\\
$C_n^{(\lambda)}(x)$ & \cite{as}, 22.3.4 & $(-1)^k  2^{n-2k}\frac{\Gamma(n - k + \lambda)}{\Gamma(\lambda) k! (n - 2k)!}$\\
&&\\
$T_0(x)$ & \cite{as}, 22.4.4 & $1$\\
$T_n(x)$ ($n > 0$) & \cite{as}, 22.3.6 &$(-1)^k  2^{n-2k} \frac{n}{2} \frac{(n - k - 1)! }{k! (n - 2k)!}$\\
&&\\
$U_n(x)$  & \cite{as}, 22.3.7 & $(-1)^k  2^{n-2k} {{n-k} \choose k}$\\
&&\\
$P_n(x)$ &  \cite{as}, 22.3.8& $(-1)^k  2^{-n} {n \choose k} {{2n - 2k} \choose n}$\\
&&\\
$H_n(x)$ & \cite{as}, 22.3.10 & $(-1)^k  2^{n-2k} \frac{n!}{k! (n - 2k)!}$\\
&&\\
$He_n(x)$ & \cite{as}, 22.3.11 & $(-1)^k  2^{-k}\frac{n!}{k! (n - 2k)!}$ \\
&&\\\hline
\end{tabular}
\caption{Mappings to the Monomials using $\beta$} 
\end{table} 

\begin{table}[h!]
\centering
\resizebox{\textwidth}{!}{
\begin{tabular}{|l|l|l|}\hline
Domain Basis & Equation & $\alpha(n, k)$\\ \hline
&&\\
$P_n^{(\alpha, \beta)}(x)$ & (\ref{Pz1})& $\frac{(-1)^{n-k}\Gamma(\alpha+n+1)}{n!\Gamma(\alpha+\beta+n+1)}$ \\
&&$\sum\limits_{l=0}^{k}(-2)^{l - n}   {{n-l} \choose {n-k}} {n \choose {n-l}}  \frac{\Gamma(\alpha+\beta+2n - l +1)}{\Gamma(\alpha + n-l + 1)}$\\
&& \\
$V_n(x)$ &  (\ref{Vx})  & $\frac{2^{n}}{{{2n} \choose n}}  \frac{(-1)^{k}}{ (n-k)!} \sum\limits_{l=0}^{k} 2^l    \frac{(1 +n)_{n-l}(\frac{1}{2} - n)_{l}}{(k-l)! \; l!}$ \\
& & \\
$W_n(x)$ &  (\ref{Wx})  & $\frac{2^{n}}{{{2n} \choose n}}  \frac{(-1)^{k}}{(n-k)!}  \sum\limits_{l=0}^{k} 2^l    \frac{(1 +n)_{n-l}(-\frac{1}{2} - n)_{l}}{(k-l)! \;  l!}$ \\
&&\\
$T_0^{\ast}(x)$ & (\ref{MTs}) & $1$\\
$T_n^{\ast}(x)$  ($n > 0$) & (\ref{MTs}) & $ \frac{n}2 \sum\limits_{v=0}^{\lfloor \frac{k}2 \rfloor}  {{n-2v} \choose {k-2v}}  (-1)^{k-v} 2^{2(n-v)-k}    \frac{(n-v-1)!}{v!(n-2v)!}$\\
&&\\
$U_n^{\ast}(x)$ & (\ref{MUs}) & $\sum\limits_{v=0}^{\lfloor \frac{k}2 \rfloor}  {{n-2v} \choose {k-2v}} {{n-v}\choose v} (-1)^{k-v} 2^{2(n-v)-k}$\\
&&\\
$V_n^{\ast}(x)$ & (\ref{MVs}) & $(-1)^{k} \frac{2^{2n-k}}{{{2n} \choose n} (n-k)!}  \sum\limits_{v=0}^k
  \sum\limits_{l=0}^{v} 2^l    \frac{(1 +n)_{n-l}(\frac{1}{2} - n)_{l}}{(k-v)! \; (v-l)! \;  l!}$\\
&&\\
$W_n^{\ast}(x)$ & (\ref{MWs}) & $(-1)^{k} \frac{2^{2n-k}}{{{2n} \choose n} (n-k)!}  \sum\limits_{v=0}^k
  \sum\limits_{l=0}^{v} 2^l    \frac{(1 +n)_{n-l}(-\frac{1}{2} - n)_{l}}{(k-v)! \; (v-l)! \;  l!}$\\
&&\\
$P_n^{\ast}(x)$ & (\ref{MPs}) & $2^{-k} \sum\limits_{v=0}^{\lfloor \frac{k}2 \rfloor}  {{n-2v} \choose {k-2v}} {{2n-2v}\choose n} {n \choose v} (-1)^{k-v} $\\
&&\\
$L_n^{(\alpha)}(x)$ & \cite{as}, 22.3.9 & $\frac{(-1)^{n-k}}{(n-k)!} {{n + \alpha} \choose k} $ \\ 
&&\\\hline
\end{tabular}}
\caption{Mappings to the Monomials using $\alpha$} 
\end{table}

\begin{table}[h!]
\centering
\begin{tabular}{|l|l|l|}\hline
Range Basis & Equation & $\beta(n, k)$\\ \hline
&&\\
$C_n^{(\lambda)}(x)$ & \cite{dlmf}, 18.18.17 & $2^{-n} \frac{\lambda + n - 2k}{\lambda} \frac{n!}{(\lambda +1)_{n-k} k!}$\\
&&\\
$T_n(x)$ & \cite{handscomb} (2.14)  & $2^{1-n} {n \choose k}$ $\mbox{ if $2 k< n$, and}$ $2^{-n} {n \choose k}$ $\mbox{ if $2k  = n$}$ \\
&&\\
$U_n(x)$  & \cite{dlmf}, 18.18.17&$2^{-n} (n - 2k + 1) \frac{n!}{(2)_{n-k} k!}$\\
&&\\
$P_n(x)$ &  \cite{dlmf}, 18.18.17& $2^{-n} (2(n - 2k) + 1) \frac{n!}{(\frac32)_{n-k} k!}$\\
&&\\
$H_n(x)$ & \cite{dlmf},  18.18.20 & $2^{-n} \frac{n!}{ (n - 2k)!k!}$\\
&&\\
$He_n(x)$ & \cite{wiki-hermite}  & $2^{-k}\frac{n!}{ (n-2k)! k!}$ \\
&&\\
\hline
\end{tabular}
\caption{Mappings from the Monomials using $\beta$} 
\end{table}

\begin{table}[h!]
\centering
\resizebox{\textwidth}{!}{
\begin{tabular}{|l|l|l|}\hline
Range Basis & Equation & $\alpha(n, k)$\\ \hline
&&\\
$P_n^{(\alpha, \beta)}(x)$ & (\ref{zP})& $(\alpha+\beta +2(n-k) + 1 )(\alpha + \beta+2)_{n-k-1}$\\
& & $\sum\limits_{l = 0}^{k} 2^{n-l} {n \choose {n-l}} {(-1)}^l  \frac{(\beta+1 + n-k)_{k-l}}{(\alpha + \beta + 2)_{n-l}} \frac{ (k - l  + 1)_{n-k}}{ (n-l+ \alpha + \beta + 2)_{n-k}}$\\
&& \\
$V_n(x)$ &  (\ref{xV})  & $n! \sum\limits_{l = 0}^{k} \frac{ (-1)^l  \;(2(n-l) +1)!!}{(k-l)! \; (2n- l -k  +1)! \; l!}$ \\
& & \\
$W_n(x)$ &  (\ref{xW})  & $(2(n-k)+1) n!\sum\limits_{l = 0}^{k}  \frac{(-1)^l \; (2(n-l)-1)!! }{(k-l)! \; (2n-l -k +1)! \;  l!}$ \\
&&\\
$T_n^{\ast}(x)$    & (\ref{Ts1M}) & $2^{1-2n} \sum\limits^{k}_{\mathclap{\substack{v = 0 \\ k-v \mbox{\small \; even}}}} {{n-v} \choose { \frac{k-v}2}} {n \choose v}  2^v$ $\mbox{ if $ k< n$} $\\
& & \\
$T_n^{\ast}(x)$ & (\ref{Ts2M}) &  $2^{-2n} \sum\limits^{k}_{\mathclap{\substack{v = 0 \\ k-v \mbox{\small \; even}}}} {{n-v} \choose { \frac{k-v}2}}    {n \choose v}  2^v $  $\mbox{ if $k  = n$}$\\
&&\\
$U_n^{\ast}(x)$ & (\ref{UsM}) & $2^{-2n} (n -k+1) n! \sum\limits^{k}_{\mathclap{\substack{v = 0 \\ k-v \mbox{\small \; even}}}} \frac{ 2^{v} }{(2)_{n - \frac{v+k}2} v! (\frac{k-v}2)!}$\\
&&\\
$V_n^{\ast}(x)$ & (\ref{VsM}) & $2^{-n} n! \sum\limits_{v=0}^k
  \sum\limits_{l = 0}^{k-v} \frac{ (-1)^l  \;(2(n -v-l) +1)!!}{(k-v- l)! \; (2n-v - l - k  +1)! \; v! \; l!}$ \\
&&\\
$W_n^{\ast}(x)$ & (\ref{WsM}) & $2^{-n} (2(n-k) + 1) n! \sum\limits_{v=0}^k 
  \sum\limits_{l = 0}^{k-v} \frac{ (-1)^l  \;(2(n -v-l) -1)!!}{(k-v- l)! \; (2n-v - l - k  +1)! \; v! \; l!}$\\
  &&\\
$P_n^{\ast}(x)$ & (\ref{PsM}) & $2^{-2n}(2(n-k)+1) n! \sum\limits^{k}_{\mathclap{\substack{v = 0 \\ k-v \mbox{\small \; even}}}}    \frac{2^{v}}{(\frac32)_{n - \frac{v+k}2} v! (\frac{k-v}2)!} $\\
&&\\
$L_n^{(\alpha)}(x)$ & \cite{dlmf}, 18.18.19 & $(-n)_{n-k} (\alpha+n-k + 1)_k $ \\ 
&&\\\hline
\end{tabular}}
\caption{Mappings from the Monomials using $\alpha$} 
\end{table}

\end{document}